\documentclass{amsart}
\usepackage[T1]{fontenc}
\usepackage[utf8]{inputenc}
\usepackage[english]{babel}
\usepackage{csquotes}

\usepackage{biblatex}
\renewbibmacro{in:}{%
        \ifentrytype{article}{}{\printtext{\bibstring{in}\intitlepunct}}}

\usepackage[hang,flushmargin]{footmisc} 

\usepackage{amsmath}
\usepackage{amssymb}
\usepackage{amsthm}

\usepackage{abstract}
\addto\captionsenglish{}

\newtheorem{defn}{Definition}[section]
\newtheorem{theo}[defn]{Theorem}

\newtheorem{lemma}[defn]{Lemma}
\newtheorem{cor}[defn]{Corollary}

\newcommand{\R}{\mathbb{R}}

\newcommand{\func}[3]{#1:#2\rightarrow#3}
\newcommand{\norm}[2]{{||#1||}_{#2}}

\newcommand{\cl}{\textup{cl}}
\newcommand{\ir}{\textup{ir}}
\newcommand{\kr}{\textup{kr}}
\newcommand{\fr}{\textup{fr}}
\newcommand{\lr}{\textup{lr}}
\newcommand{\bir}{\textup{bir}}
\newcommand{\bkr}{\textup{bkr}}
\newcommand{\bfr}{\textup{bfr}}
\newcommand{\blr}{\textup{blr}}
\newcommand{\wscl}{\textup{w}^*\textup{-cl}}

\renewcommand{\epsilon}{\varepsilon}
\renewcommand{\theta}{\vartheta}
\renewcommand{\phi}{\varphi}

\title{A Topological Characterisation of Haar Null Convex Sets}
\author{Davide Ravasini}
\address{Institut für Mathematik, Universität Innsbruck \newline\indent Technikerstraße 13, 6020 Innsbruck, Austria}
\email{davide.ravasini@uibk.ac.at}

\begin{document}
\maketitle
\let\thefootnote\relax\footnote{\today \newline \indent\emph{2020 Mathematics Subject Classification:} 46B20 (primary), 52A07 (secondary). \newline
\indent \emph{Keywords:} Haar null set, loose radius, weak$^*$ closure. \newline 
\indent The author's work is supported by the Austrian Science Fund (FWF): P 32523-N.}

\begin{abstract}
\noindent \textsc{Abstract}. In $\R^d$, a closed, convex set has zero Lebesgue measure if and only its interior is empty. More generally, in separable, reflexive Banach spaces, closed and convex sets are Haar null if and only if their interior is empty. We extend this facts by showing that a closed, convex set in a separable Banach space is Haar null if and only if its weak$^*$ closure in the second dual has empty interior with respect to the norm topology. It then follows that, in the metric space of all nonempty, closed, convex and bounded subsets of a separable Banach space, converging sequences of Haar null sets have Haar null limits.
\end{abstract}

\section{Introduction}
A Borel set $E$ in an Abelian Polish group $X$ is said to be \emph{Haar Null} if there is a Borel probability measure $\mu$ on $X$ such that $\mu(x+E)=0$ for every $x\in X$. Haar null sets were introduced for the first time by J.P.R.\ Christensen in \cite{chris_1972} in order to extend the notion of sets with zero Haar measure to nonlocally compact Polish groups, where the Haar measure is not defined. Indeed, Haar null sets and sets with zero Haar measure agree on locally compact Abelian Polish groups and, as in the locally compact case, Haar null sets form a $\sigma$-ideal in the Borel $\sigma$-algebra of $X$. We say that a Borel set $E\subseteq X$ is \emph{Haar positive} if $E$ is not Haar null. The reader is invited to have a look at the survey papers \cite{bog_2018} and \cite{en_2020}, as well as at \cite{benlin}, Chapter 6, for a detailed exposition of this topic.

Although the definition of Haar positive sets is measure-theoretical in nature, in \cite{daverava2} a measure-free characterisation of Haar positive closed, convex sets in separable Banach spaces is provided.

\begin{theo}
\label{theo:intro}
Let $C$ be a closed and convex set in a separable Banach space $X$ with unit ball $B_X$. The following assertions are equivalent.
\begin{enumerate}
\item $C$ is Haar positive.
\item There is $r>0$ with the property that, for every compact set $K\subseteq rB_X$, there is $x\in X$ such that $x+K\subseteq C$.
\end{enumerate}
\end{theo} 
\noindent Equivalently, a closed, convex set $C$ in a separable Banach space $X$ is Haar null if and only if $C$ is \emph{Haar meagre}. That is, there is a compact metric space $M$ and a continuous function $\func{f}{M}{X}$ such that $f^{-1}(x+C)$ is meagre in $M$ for every $x\in X$. For a more general treatment of Haar meagre sets, we refer the reader to \cite{dar_2013}, where they were first introduced.

Theorem \ref{theo:intro} motivates the introduction of several geometric radii associated to such sets. These are defined in Section \ref{sec:radii} and multiple properties about these quantities are shown. In Section \ref{sec:weaks} we exploit these radii to prove a new characterisation of Haar positive, closed, convex sets. Namely, a closed, convex set $C$ in a separable Banach space $X$ is Haar positive if and only if its weak$^*$ closure in the second dual $X^{**}$ has nonempty interior with respect to the norm topology. This improves the well-known fact that closed, convex subsets of a Euclidean space $\R^d$ have positive Lebesgue measure if and only if their interior is nonempty and a theorem of Eva Matou\v{s}kov\'{a} (\cite{eva3}) which states that, in separable, reflexive Banach spaces, closed convex sets are Haar positive if and only if they have nonempty interior. As a corollary, it is shown in Section \ref{sec:top} that the family of Haar positive, closed, convex and bounded sets is open in the space of all nonempty, closed, convex and bounded subsets of $X$, endowed with the Hausdorff distance.

The standard notation of Banach space theory is used throughout the paper. Given a Banach space $X$, $B_X$ and $S_X$ stand for the closed unit ball and the unit sphere of $X$ respectively. $X^*$ denotes the dual of $X$, whereas $X^{**}$ is the second dual. The closure of a set $E\subseteq X$ is denoted by $\cl(E)$ and in a dual space we denote by $\wscl(E)$ the closure of $E$ in the weak$^*$ topology. We use the notation $\mathcal{C}(X)$ for the space of all nonempty, closed, convex and bounded subsets of $X$. This turns into a complete metric space if endowed with the Hausdorff distance $d_\textup{H}$ given by
\[ d_\textup{H}(C,D)=\inf\{\epsilon>0\,:\,C\subseteq\cl(D+\epsilon B_X) \text{ and }D\subseteq\cl(C+\epsilon B_X)\}. \]
All Banach spaces are assumed to be real.

\section{Geometric radii of closed, convex and bounded sets} 
\label{sec:radii}
Given $\rho>0$ and a nonempty, closed, convex set $C$ in a Banach space $X$, we denote by $\ir(C,\rho)$ the $\rho\,$-\emph{inner radius} of $C$. That is, the supremum of all $r\geq 0$ such that $C$ contains a closed ball of the form $x+rB_X$, where $x\in\rho B_X$. Clearly, $C$ has nonempty interior if and only if $\ir(C,\rho)>0$ for some $\rho>0$. The $\rho\,$-\emph{compact radius} $\kr(C,\rho)$ is defined as follows: it is the supremum of all $r\geq 0$ with the property that, for every compact set $K\subseteq rB_X$, there is $x\in\rho B_X$ such that $x+K\subseteq C$. The $\rho\,$-\emph{finite radius} $\fr(C,\rho)$ is defined similarly. Namely, it is the supremum of all $r\geq 0$ with the property that, for every finite set $F\subseteq rB_X$, there is $x\in\rho B_X$ such that $x+F\subseteq C$. Finally, we introduce the $\rho\,$-\emph{loose radius} $\lr(C,\rho)$ as follows: it is the supremum of all $r\geq 0$ with the property that, for every finite set $F=\{x_1,\dots,x_n\}\subseteq rB_X$ and every $\epsilon>0$, there are a finite set $G=\{y_1,\dots,y_n\}$ and $z\in\rho B_X$ such that $\norm{x_j-y_j}{X}<\epsilon$ for every $j\in\{1,\dots,n\}$ and $z+G\subseteq C$.

\begin{theo}
\label{theo:radii}
Let $C$ be a nonempty, closed, convex set in a Banach space $X$. The following inequalities hold for every $\rho>0$.
\[ \ir(C,\rho)\leq\kr(C,\rho)\leq\fr(C,\rho)\leq\lr(C,\rho)\leq 2\kr(C,\rho). \]
\end{theo}
\begin{proof}
The only nontrivial inequality is the last one and we will prove it using a variation of the Banach-Dieudonn\'{e} Theorem, similar to the one which appears in \cite{eva1}. If $\lr(C,\rho)=0$, the claim is obvious. Assume that $\lr(C,\rho)>0$. We aim to show that, for every $r\geq 0$ such that $2r<\lr(C,\rho)$, every compact set $K\subseteq rB_X$ can be translated into $C$ via some $z\in\rho B_X$. Let $K$ be a compact subset of $rB_X$ and find a finite set $F_1=\{x_{1,1},\dots,x_{1,k(1)}\}\subseteq 2rB_X$ such that $2^{-1}F_1$ is an $(r/4)$-net for $K$. Find a further finite set $G_1=\{y_{1,1},\dots,y_{1,k(1)}\}$ and $z_1\in\rho B_X$ such that $\norm{y_{1,j}-x_{1,j}}{X}<r/2$ for every $j\in\{1,\dots,k(1)\}$ and such that $z_1+G_1\subseteq C$. Notice that $2^{-1}G_1$ is an $(r/2)$-net of $K$. Define 
\[ K_1=2\Bigl((K-2^{-1}G_1)\cap\frac{r}{2}B_X\Bigr), \]
and observe that $K_1$ is a compact subset of $rB_X$ such that $K\subseteq 2^{-1}G_1+2^{-1}K_1$. Working by induction, let $n$ be a positive integer and assume that we have already found finite sets $G_1,\dots,G_n\subseteq X$, $z_1,\dots,z_n\in\rho B_X$ and a compact set $K_n\subset rB_X$ such that $z_j+G_j\subseteq C$ for every $j\in\{1,\dots,n\}$ and 
\begin{equation}
\label{eq:loose1}
K\subseteq \sum_{j=1}^n2^{-j}G_j+2^{-n}K_n.
\end{equation}
Let $F_{n+1}=\{x_{n+1,1},\dots,x_{n+1,k(n+1)}\}\subseteq 2rB_X$ be a finite set such that $2^{-1}F_{n+1}$ is an $(r/4)$-net for $K_n$. Find a further finite set $G_{n+1}=\{y_{n+1,1},\dots,y_{n+1,k(n+1)}\}$ and $z_{n+1}\in\rho B_X$ such that $\norm{y_{n+1,j}-x_{n+1,j}}{X}<r/2$ for every $j\in\{1,\dots,k(n+1)\}$ and $z_{n+1}+G_{n+1}\subseteq C$. Notice again that $2^{-1}G_{n+1}$ is an $(r/2)$-net of $K_n$. Define 
\[ K_{n+1}=2\Bigl((K_n-2^{-1}G_{n+1})\cap\frac{r}{2}B_X\Bigr), \]
and observe that $K_{n+1}\subseteq rB_X$ is compact and that $K_n\subseteq 2^{-1}G_{n+1}+2^{-1}K_{n+1}$ holds. Now we have 
\[ K\subseteq \sum_{j=1}^{n+1}2^{-j}G_j+2^{-(n+1)}K_{n+1}, \]
therefore we can claim by induction that we can find a sequence ${(z_n)}_{n=1}^\infty$ in $\rho B_X$, a sequence ${(G_n)}_{n=1}^\infty$ of finite subsets of $X$ and a sequence ${(K_n)}_{n=1}^\infty$ of compact subsets of $rB_X$ such that $z_n+G_n\subseteq C$ and (\ref{eq:loose1}) hold for every $n$. In particular, this implies 
\[ K\subseteq\cl\biggl(\,\sum_{n=1}^\infty 2^{-n}G_n\biggr). \]
Put $z=\sum_{n=1}^\infty 2^{-n}z_n$ and observe that $z\in\rho B_X$. Then
\[ z+K\subseteq\sum_{n=1}^\infty 2^{-n}z_n+\cl\biggl(\,\sum_{n=1}^\infty 2^{-n}G_n\biggr)=\cl\biggl(\,\sum_{n=1}^\infty 2^{-n}z_n+\sum_{n=1}^\infty 2^{-n}G_n\biggr)\subseteq C, \]
where the last inclusion holds because $C$ is closed and convex. Since $K\subseteq rB_X$ is arbitrary, we are done.
\end{proof}

The following argument also appears with some minor differences in \cite{eva1}. We report it here for completeness and formulate it in the language we have just introduced.
\begin{lemma}
\label{lemma:kr}
Let $C$ be a closed, convex set in a separable Banach space $X$. $C$ is Haar positive if and only if there is $\rho>0$ such that $\kr(C,\rho)>0$. It follows in particular that $C$ is Haar positive if and only if $\lr(C,\rho)>0$ for some $\rho>0$.
\end{lemma}
\begin{proof}
It is a direct consequence of Theorem \ref{theo:intro} that $\kr(C,\rho)=0$ for every $\rho>0$ if $C$ is Haar null. To show the converse statement, assume by way of contradiction that, for every positive integer $n$, there is a compact set $K_n\subseteq n^{-1}B_X$ such that, whenever $x+K_n\subseteq C$ for some $x\in X$, we have $\norm{x}{X}>n$. By Theorem \ref{theo:intro}, there is $r>0$ such that every compact subset of $rB_X$ can be translated into $C$. In particular, this holds for
\[ K=\{0\}\cup\bigcup_{n>r^{-1}}K_n. \]
If $x$ is such that $x+K\subseteq C$, then $x+K_n\subseteq C$ for every $n>r^{-1}$. It then follows by the choice of $K_n$ that $\norm{x}{X}>n$ for every positive integer $n$, which does not make any sense. The last assertion is a consequence of Theorem \ref{theo:radii}.
\end{proof}

\begin{lemma}
\label{lemma:lr-fr}
Given $\rho>0$ and a nonempty, closed, convex set $C$ in a Banach space $X$, we have
\[ \lr(C,\rho)=\lim_{\delta\to 0^+}\fr\bigl(\cl(C+\delta B_X),\rho\bigr). \]
\end{lemma}
\begin{proof}
Set $r=\lr(C,\rho)$. Pick $\epsilon>0$ and find $F=\{x_1,\dots,x_n\}\subset(r+\epsilon)B_X$ and $\delta_0>0$ such that every finite set $G=\{y_1,\dots,y_n\}$ which fulfills the condition $\norm{y_j-x_j}{X}<2\delta_0$ for every $j\in\{1,\dots,n\}$ cannot be translated into $C$ via some $z\in\rho B_X$. Suppose that there is $z\in\rho B_X$ such that $z+F\subset\cl(C+\delta_0B_X)$. This would imply that, for every $j\in\{1,\dots,n\}$, we can find $w_j\in C$ such that $\norm{z+x_j-w_j}{X}<2\delta_0$. Put $y_j=w_j-z$ for every $j$. Then the set $G=\{y_1,\dots,y_n\}$ satisfies $z+G\subseteq C$ and $\norm{x_j-y_j}{X}<2\delta_0$ for every $j$, in contradiction with the choice of $F$ and $\delta_0$. Thus, $F$ cannot be translated into $\cl(C+\delta_0B_X)$ via some $z\in\rho B_X$, which yields
\[ \lim_{\delta\to 0^+}\fr\bigl(\cl(C+\delta B_X),\rho\bigr)\leq\fr\bigl(\cl(C+\delta_0B_X),\rho\bigr)\leq r+\epsilon. \]
Since $\epsilon$ is arbitrary, we get
\[ \lim_{\delta\to 0^+}\fr\bigl(\cl(C+\delta B_X),\rho\bigr)\leq r. \]
To show the opposite inequality, observe that it is obvious if $r=0$. Assume that $r>0$, choose $\delta\in(0,r)$ and pick a finite set $F=\{x_1,\dots,x_n\}\subset(r-\delta)B_X$. Find $G=\{y_1,\dots,y_n\}$ and $z\in\rho B_X$ such that $z+G\subseteq C$ and $\norm{x_j-y_j}{X}<\delta$ for each $j$. Notice that
\[ z+x_j=z+y_j+x_j-y_j\in C+\delta B_X. \]
for each $j$, hence $z+F\subseteq\cl(C+\delta B_X)$. Since $\delta>0$ and $F\subset (r-\delta)B_X$ are arbitrary, we conclude that
\[ \lim_{\delta\to 0^+}\fr\bigl(\cl(C+\delta B_X),\rho\bigr)\geq\lim_{\delta\to 0^+}(r-\delta)=r, \]
as wished.
\end{proof}

\section{Weak$^*$ closures of Haar positive closed, convex sets}
\label{sec:weaks}
Given a Banach space $X$ and a positive integer $n$, recall that we can endow the Banach space $X^n$, the product of $n$ copies of $X$, with the $\infty$-product norm:
\[ \norm{x}{X^n}=\max_{1\leq j\leq n}\norm{x_j}{X} \]
for every $x=(x_1,\dots,x_n)\in X^n$. In this way we have $B_{X^n}={(B_X)}^n$. We denote by $\func{\Delta}{X}{X^n}$ the diagonal embedding $x\mapsto(x,\dots,x)$. The second dual space of $X^n$ is simply ${(X^{**})}^n$, endowed with the same $\infty$-product norm. Notice that, in ${(X^{**})}^n$, we have $\wscl(\Delta(\rho B_X))=\Delta(\rho B_{X^{**}})$. We are now ready to state our main theorem.
\begin{theo}
\label{theo:hpws}
Let $\rho>0$, let $C$ be a nonempty, closed, convex set in a Banach space $X$ and let $r\geq 0$. The following assertions are equivalent.
\begin{enumerate}
\item $\lr(C,\rho)\geq r$.
\item For every finite $F\subseteq rB_X$ there is $z\in\rho B_{X^{**}}$ such that $z+F\subseteq\wscl(C)$.
\item There is $z_0\in\rho B_{X^{**}}$ such that $z_0+rB_{X^{**}}\subseteq\wscl(C)$.
\end{enumerate}
In particular, $\lr(C,\rho)=\ir(\wscl(C),\rho)$ and, in case $X$ is separable, $C$ is Haar positive if and only if $\wscl(C)$ has nonempty interior in the norm topology of $X^{**}$.
\end{theo}
\begin{proof}
(1)$\Rightarrow$(2). Let $F=\{x_1,\dots,x_n\}\subseteq rB_X$ be a finite set and, for each positive integer $k$, find $G_k=\{y_{k,1},\dots,y_{k,n}\}$ and $z_k\in\rho B_X$ such that $\norm{x_j-y_{k,j}}{X}<k^{-1}$ for each $j\in\{1,\dots,n\}$ and $z_k+G_k\subseteq C$. Since ${(z_k)}_{k=1}^\infty$ is a bounded sequence, it admits a subnet with a weak$^*$ limit $z\in\rho B_{X^{**}}$. Now, it is not hard to see that $z+F\subseteq\wscl(C)$.

(2)$\Rightarrow$(3). Let $\mathcal{D}$ be the directed set consisting of all finite subsets of $rB_X$, ordered by inclusion. For every $F\in\mathcal{D}$, let $z_F\in\rho B_{X^{**}}$ be such that $z_F+F\subseteq\wscl(C)$. Since ${(z_F)}_{F\in\mathcal{D}}$ is a net in the w$^*$-compact set $\rho B_{X^{**}}$, there is a directed set $I$ and a map $\func{\psi}{I}{\mathcal{D}}$ such that ${(z_{\psi(\alpha)})}_{\alpha\in I}$ is a subnet of ${(z_F)}_{F\in\mathcal{D}}$ with a weak$^*$ limit $z_0\in\rho B_{X^{**}}$. We want to show that $z_0+rB_{X^{**}}$ is contained in $\wscl(C)$. It suffices to show that $z_0+rB_X$ is contained in $\wscl(C)$, because $rB_{X^{**}}$ is the weak$^*$ closure of $rB_X$. Given $x\in rB_X$, the set $J=\{\alpha\in I\,:\,x\in\psi(\alpha)\}$ is a directed subset of $I$. The net ${(z_{\psi(\alpha)})}_{\alpha\in J}$ is a subnet of ${(z_{\psi(\alpha)})}_{\alpha\in I}$ and has therefore the same weak$^*$ limit $z_0$. Since $z_{\psi(\alpha)}+\psi(\alpha)\subseteq\wscl(C)$ for every $\alpha\in J$, we have in particular $z_{\psi(\alpha)}+x\in\wscl(C)$ for every $\alpha\in J$. By taking the limit, we get $z_0+x\in\wscl(C)$. Since $x$ is arbitrary, the inclusion $z_0+rB_X\subseteq\wscl(C)$ is proved.

(3)$\Rightarrow$(1). Assume by contradiction that $\lr(C,\rho)<r$. By Lemma \ref{lemma:lr-fr}, there is $\delta>0$ such that $\fr(\cl(C+\delta B_X),\rho)<r$, which in turn means that there is a finite set $F=\{x_1,\dots,x_n\}\subseteq rB_X$ such that $(z+F)\setminus\cl(C+\delta B_X)\ne\varnothing$ for every $z\in\rho B_X$. In $X^n$ define $D=(C-x_1)\times\cdots\times(C-x_n)$ and observe that
\[ \cl(D+\delta B_{X^n})=\bigl(\cl(C+\delta B_X)-x_1\bigr)\times\cdots\times\bigl(\cl(C+\delta B_X)-x_n\bigr). \]
Notice moreover that $\Delta(\rho B_X)\cap\cl(D+\delta B_{X^n})=\varnothing$, otherwise it would hold that $z+F\subseteq\cl(C+\delta B_X)$ for every $(z,\dots,z)\in \Delta(\rho B_X)\cap\cl(D+\delta B_{X^n})$, against how $F$ has been chosen. Hence, by the Hahn-Banach theorem, there are $f\in{(X^n)}^*$ and $t_1,t_2\in\R$ such that $t_1<t_2$, $f(x)<t_1$ for every $x\in\Delta(\rho B_X)$ and $f(x)>t_2$ for every $x\in D$. By taking the weak$^*$ closures of both sets in ${(X^n)}^{**}$, we have $f(x)\leq t_1$ for every $x\in\Delta(\rho B_{X^{**}})$ and $f(x)\geq t_2$ for every $x\in\wscl(D)$, thus $\Delta(\rho B_{X^{**}})\cap\wscl(D)=\varnothing$. This is a contradiction, because
\[ z_0\in\bigcap_{j=1}^n\bigl(\wscl(C)-x_j\bigr), \]
i.e.\ $(z_0,\dots,z_0)\in(\wscl(C)-x_1)\times\cdots\times(\wscl(C)-x_n)=\wscl(D)$.

The last statement is a consequence of Lemma \ref{lemma:kr}.
\end{proof}

Theorem \ref{theo:hpws} offers an interesting connection with the theory of weak$^*$ derived sets. If $X$ is a Banach space and $E$ is a subset of $X^*$, the first weak$^*$ derived set of $E$ is given by
\[ E^{(1)}=\bigcup_{n=1}^\infty\wscl(E\cap nB_{X^*}) \]
and corresponds to the set of all possible limits of bounded, weak$^*$ convergent nets with elements in $E$. Clearly $E^{(1)}\subseteq\wscl(E)$. We refer the reader to \cite{ostr_2023} and the references therein for a detailed account on weak$^*$ derived sets. In particular, it has to be remarked that, in general, $E^{(1)}$ and $\wscl(E)$ can be different. However, we have the following result.

\begin{cor}
Let $C$ be a closed, convex set in a Banach space $X$. In the second dual $X^{**}$, $C^{(1)}$ has empty interior in the norm topology if and only if $\wscl(C)$ does.
\end{cor}
\begin{proof}
One of the implications follows from $C^{(1)}\subseteq\wscl(C)$. Conversely, assume that $\wscl(C)$ has nonempty interior in the norm topology. Then, by Theorem \ref{theo:hpws} and Theorem \ref{theo:radii}, there are $\rho>0$ and $r>0$ such that $\kr(C,\rho)\geq r$. Now, it is not hard to see that $\kr(C\cap nB_X,\rho)\geq r$ for every integer $n>\rho+r$, hence $\wscl(C\cap n B_X)$ has nonempty interior in the norm topology for every $n>\rho+r$ and so does $C^{(1)}$.
\end{proof}

\section{Haar positive sets and the Hausdorff metric}
\label{sec:top}

Let $C$ be a closed and convex set in a Banach space $X$. Under the additional assumption that $C$ is bounded we define 
\begin{align*} 
\bir(C)=\sup_{\rho>0}\ir(C,\rho),&\quad\bkr(C)=\sup_{\rho>0}\kr(C,\rho),\\
\bfr(C)=\sup_{\rho>0}\fr(C,\rho),&\quad\blr(C)=\sup_{\rho>0}\lr(C,\rho). 
\end{align*}
All these values are finite, as $\textup{diam}(C)$ is an upper estimate for all of them. Moreover, the chain of inequalities 
\begin{equation} 
\label{eq:radii}
\bir(C)\leq\bkr(C)\leq\bfr(C)\leq\blr(C)\leq 2\bkr(C)
\end{equation}
is a direct consequence of Theorem \ref{theo:radii}. We call $\bir(C)$ the \emph{bounded inner radius} of $C$. It is the supremum of all $r\geq 0$ such that $C$ contains a closed ball of radius $r$. $\bkr(C)$ is the \emph{bounded compact radius} of $C$ and corresponds to the supremum of all $r\geq 0$ with the property that, for every compact set $K\subseteq rB_X$, there is $x\in X$ such that $x+K\subseteq C$. The \emph{bounded finite radius} $\bfr(C)$ is, similarly, the supremum of all $r\geq 0$ with the property that, for every finite set $F\subseteq rB_X$, there is $x\in X$ such that $x+F\subseteq C$. Finally, the \emph{bounded loose radius} $\blr(C)$ is the supremum of all $r\geq 0$ with the property that, for every finite set $F=\{x_1,\dots,x_n\}\subseteq rB_X$ and every $\epsilon>0$, there are a finite set $G=\{y_1,\dots,y_n\}$ and $z\in X$ such that $\norm{x_j-y_j}{X}<\epsilon$ for every $j$ and $z+G\subseteq C$. In case $X$ is separable, $C$ is Haar positive if and only if $\blr(C)>0$. This follows from ($\ref{eq:radii}$) and Theorem \ref{theo:intro}.

Theorem \ref{theo:hpws} allows to prove that the family $\mathcal{H}^+(X)$ of Haar positive, closed, convex and bounded subsets of a separable Banach space $X$ is an open subset of $\mathcal{C}(X)$. This result follows from a few lemmas we are going to show.

\begin{lemma}
\label{lemma:hpws}
Given a closed, convex and bounded set $C$ in a Banach space $X$, we have $\blr(C)=\bir(\wscl(C))$.
\end{lemma}
\begin{proof}
Using Theorem \ref{theo:hpws}, we have
\[ \blr(C)=\sup_{\rho>0}\lr(C,\rho)=\sup_{\rho>0}\ir\bigl(\wscl(C),\rho\bigl)=\bir\bigl(\wscl(C)\bigr). \qedhere \]
\end{proof}

\begin{lemma}
\label{lemma:wscliso}
Given a Banach space $X$, the function $\func{\wscl}{\mathcal{C}(X)}{\mathcal{C}(X^{**})}$ is an isometry.
\end{lemma}
\begin{proof}
Take $\epsilon>0$ and let $C,D\in\mathcal{C}(X)$ be such that $d_\textup{H}(C,D)\leq\epsilon$. This means that $C\subseteq\cl(D+\epsilon B_X)$ and $D\subseteq\cl(C+\epsilon B_X)$. Take $x\in\wscl(C)$ and let ${(x_\alpha)}_{\alpha\in I}$ be a net in $C$ whose weak$^*$ limit is $x$. Choose $\delta>0$ and, for each $\alpha\in I$, find $y_\alpha\in D$ and $z_\alpha\in\epsilon B_X$ such that $\norm{x_\alpha-(y_\alpha+z_\alpha)}{X}\leq\delta$. By considering a subnet if necessary, we can assume that the nets ${(y_\alpha)}_{\alpha\in I}$ and ${(z_\alpha)}_{\alpha\in I}$ have weak$^*$ limits $y\in\wscl(D)$ and $z\in\epsilon B_{X^{**}}$ respectively. Moreover, 
\[ \norm{x-(y+z)}{X^{**}}\leq\liminf_{\alpha\in I}\norm{x_\alpha-(y_\alpha+z_\alpha)}{X}\leq\delta, \]
i.e.\ $x\in\wscl(D)+(\epsilon+\delta)B_{X^{**}}$. Since $\delta$ is arbitrary, we have
\[ x\in\bigcap_{\delta>0}\bigl(\wscl(D)+(\epsilon+\delta)B_{X^{**}}\bigr)=\cl\bigl(\wscl(D)+\epsilon B_{X^{**}}\bigr). \]
As $x$ is also arbitrary, we conclude that $\wscl(C)\subseteq\cl(\wscl(D)+\epsilon B_{X^{**}})$. The inclusion $\wscl(D)\subseteq\cl(\wscl(C)+\epsilon B_{X^{**}})$ is shown symmetrically, hence we get that $d_\textup{H}(\wscl(C),\wscl(D))\leq\epsilon$. 

Conversely, suppose that $d_\textup{H}(C,D)>\epsilon$ and, by swapping $C$ and $D$ if necessary, assume that there is $x_0\in C\setminus\cl(D+\epsilon B_X)$. A standard Hahn-Banach argument shows that $x_0$ does not belong to $\wscl(D+\epsilon B_X)$ either. Now observe that 
\begin{align*} 
\cl(\wscl(D)+\epsilon B_{X^{**}}) &= \wscl(D)+\epsilon B_{X^{**}}= \\
&=\wscl(D)+\wscl(\epsilon B_{X})\subseteq\wscl(D+\epsilon B_X). 
\end{align*}
Thus we get $x_0\notin\cl(\wscl(D)+\epsilon B_{X^{**}})$, from which $d_\textup{H}(\wscl(C),\wscl(D))>\epsilon$ follows. Since $C$, $D$ and $\epsilon$ are arbitrary, the proof is complete.
\end{proof}

Finally, we want to show the continuity of the bounded inner radius in the metric space of all nonempty, closed, convex and bounded subsets of a Banach space, endowed with the Hausdorff metric. Although it seems that this result cannot be found in the literature, it might be well-know and belong to the folklore. We prove it here for the sake of completeness.
\begin{lemma}
\label{lemma:ir}
In a Banach space $X$, the function $\func{\bir}{\mathcal{C}(X)}{[0,+\infty)}$ is Lipschitz continuous with Lipschitz constant $1$.
\end{lemma}
\begin{proof}
The proof is based on the following claim: if $C$ is a nonempty, closed, convex and bounded set in $X$ and $\epsilon>0$, then
\begin{equation}
\label{eq:ir+eps}
\bir\bigl(\cl(C+\epsilon B_X)\bigr)=\bir(C)+\epsilon.
\end{equation}

Let us see first that $\bir(\cl(C+\epsilon B_X))\leq\bir(C)+\epsilon$. Set $r=\bir(C)$ and, looking for a contradiction, assume that $\bir(\cl(C+\epsilon B_X))>r+\epsilon$. Then there is $\delta>0$ such that $r+\epsilon+\delta<\bir(\cl(C+\epsilon B_X))$ and, without loss of generality, we may assume that $(r+\epsilon+\delta)B_X\subseteq\cl(C+\epsilon B_X)$. Since $(r+\delta)B_X\setminus C\ne\varnothing$, the Hahn-Banach Theorem provides $x_0\in(r+\delta)B_X$, $t\in\R$ and $f\in S_{X^*}$ such that $f(x_0)>t>f(x)$ for every $x\in C$. Further, we have $t<\norm{f}{X^*}\norm{x_0}{X}\leq r+\delta$. Find $\delta'$ such that $r+\delta-\delta'-t>0$ and $x_1\in(r+\epsilon+\delta)B_X$ such that $f(x_1)>r+\epsilon+\delta-\delta'$. Since $(r+\epsilon+\delta)B_X\subseteq\cl(C+\epsilon B_X)$, there exist $x\in C$ and $y\in\epsilon B_X$ such that
\[ \norm{x_1-(x+y)}{X}<r+\delta-\delta'-t \]
On the other hand, we have
\[ \norm{x_1-(x+y)}{X}\geq f(x_1)-f(x)-f(y)>r+\epsilon+\delta-\delta'-t-\epsilon=r+\delta-\delta'-t, \]
a contradiction. 

To see $\bir(C)+\epsilon\leq\bir(\cl(C+\epsilon B_X))$, set again $r=\bir(C)$, choose $\delta\in(0,r)$ and find $x_0\in X$ with $x_0+(r-\delta)B_X\subseteq C$. Then $x_0+(r+\epsilon-\delta)B_X\subseteq\cl(C+\epsilon B_X)$ follows. Since $\delta$ is arbitrary, it follows that $\bir(\cl(C+\epsilon B_X))\geq r+\epsilon$, as wished.

To prove the statement, pick $C,D\in\mathcal{C}(X)$ and $\epsilon>0$. If $d_\textup{H}(C,D)\leq\epsilon$, then $D\subseteq\cl(C+\epsilon B_X)$, which implies by (\ref{eq:ir+eps}) that $\bir(D)\leq\bir(C)+\epsilon$. The inequality $\bir(C)\leq\bir(D)+\epsilon$ follows similarly. Thus, $|\bir(C)-\bir(D)|\leq\epsilon$. Since $C,D$ and $\epsilon$ are arbitrary, this lets us conclude that $\bir$ is $1$-Lipschitz.
\end{proof}

\begin{theo}
Let $X$ be a Banach space. The function $\func{\blr}{\mathcal{C}(X)}{\R}$ is $1$-Lipschitz. In particular, in case $X$ is separable, the family $\mathcal{H}^+(X)$ of all Haar positive closed, convex and bounded subsets of $X$ is open. Equivalently, a convergent sequence ${(C_n)}_{n=1}^\infty\subset\mathcal{C}(X)$ of Haar null sets has a Haar null limit.
\end{theo}
\begin{proof}
We have $\blr=\bir\circ\wscl$ by Lemma \ref{lemma:hpws}, hence, by Lemma \ref{lemma:wscliso} and Lemma \ref{lemma:ir}, $\blr$ is a composition of $1$-Lipschitz maps and therefore it is $1$-Lipschitz. Since
\[ \mathcal{H}^+(X)=\blr^{-1}\bigl((0,+\infty)\bigr), \]
it follows immediately that $\mathcal{H}^+(X)$ is open.
\end{proof}

\section*{Acknowledgements}
The author would like to thank Professors Eva Kopeck\'{a} and Christian Bargetz for the many helpful remarks, and the referee for the careful reading of the manuscript and the improvements in the arguments. The author is also thankful to Professor Mikhail Ostrovskii for pointing out a mistake in an earlier version of the preprint. 

\printbibliography

\end{document}